\documentclass[12pt]{article}
\usepackage{amsmath}

\usepackage{amssymb,amsmath,latexsym,amsthm}
\usepackage[english]{babel}

\newtheorem{theorem}{Theorem}[section]
\newtheorem{lemma}{Lemma}[section]

\theoremstyle{definition}

\title{Coefficients of cyclotomic polynomials}
\author{Pingzhi Yuan \\
{\small School of Mathematics, South China Normal University
  , Guangzhou 510631, P.R.CHINA}\\
 {\small e-mail  mcsypz@mail.sysu.edu.cn}}

\date{}

\begin{document}
\baselineskip=17pt \maketitle

\section*{\center{Abstract}} Let $a(n, k)$ be the $k$-th coefficient of the $n$-th cyclotomic polynomial.
Recently, Ji,  Li and  Moree \cite{JLM09} proved that  for any
integer $m\ge1$, $\{a(mn, k)| n, k\in\mathbb{N}\}=\mathbb{Z}$. In
this paper, we improve this result and prove that for any integers
 $s>t\ge0$,
 $$\{a(ns+t, k)| n, k\in\mathbb{N}\}=\mathbb{Z}.$$
\\

\noindent 2000 Mathematics Subject Classification:11B83; 11C08

Keywords: Cyclotomic polynomials; Dirichlet's theorem; Squarefree
integers
\section{Introduction}
Let $\Phi_n(x)=\sum_{k=0}^{\varphi(n)}a(n, k)x^k$ be the $n$th
cyclotomic polynomial. The  Taylor series of $1/\Phi_n(x)$ around
$x=0$ is given by $1/\Phi_n(x)=\sum_{k=0}^{\varphi(n)}c(n, k)x^k$.
It is not difficult to show that $a(n, k)$ and $c(n, k)$ are all
integers. The coefficients $a(n,k)$ and $c(n,k)$ are quite small in
absolute value, for example for $n<105$ it is well-known that
$|a(n,k)|\le1$ and for $n<561$ we have $|c(n,k)|\le1$(see
\cite{Mo09}).  Migotti \cite{Mi83}  showed that all $a(pq, i)\in\{0,
\pm1\}$, where  $p$ and $q$ are distinct primes.  Beiter \cite{Be1}
and \cite{Be2} gave a criterion on $i$ for $a(pq, i)$ to be $0, 1$
or -1, see also Lam and Leung \cite{LL96}. Also Carlitz \cite{Ca66}
computed the number of non-zero $a(pq, i)$'s. For more information
on this topic, we refer to the beautiful survey paper of Thangadurai
\cite{Th00}. Bachman \cite{Ba1, Ba2} proved the existence of an
infinite family of $n=pqr$ with all $a(pqr, i)\in\{0, \pm1\}$, where
$p, q, r$ are distinct odd primes.

Let $m\ge1$ be a  integer. Put
$$S(m)=\{a(mn,k)|n\ge1,k\ge0\}\quad\mbox{ and}\quad
R(m)=\{c(mn,k)|n\ge1,k\ge0\}.$$

Schur poved in 1931 (in a letter to E. Landau) that $S(1)$ is not a
finite set, see Lenstra \cite{Len78}. In 1987 Suzuki \cite{Su87}
proved that $S(1)=\mathbb{Z}$. Recently, Ji, Li and Moree
\cite{JLM09}, \cite{JiL08} proved that with  $S(m)=R(m)=\mathbb{Z}$
for any integer $m\ge1$.

Let $m\ge1, s>t\ge0$ be positive  integers with $\gcd(s, t)=1$. Put
$$S(m; s, t)=\{a(m(sn+t),k)|n\ge1,k\ge0\}\quad\mbox{ and}\quad
R(m; s, t)=\{c(m(sn+t),k)|n\ge1,k\ge0\}.$$

In this note, by a slight modification of the proof in \cite{JLM09},
we prove the following generalization of the result in \cite{JLM09}.

\begin{theorem}Let $m\ge1, s>t\ge0$ be positive  integers with $\gcd(s,
t)=1$. Then $S(m; r, t)=R(m; s, t)=\mathbb{Z}$.\end{theorem}

An equivalent statement of Theorem 1.1 is the following result,
which is the motivation to write this paper.
\begin{theorem}Let $s>t\ge0$ be integers, then
$$\{a(ns+t, k)| n, k\in\mathbb{N}\}=\{c(ns+t, k)| n,
k\in\mathbb{N}\}=\mathbb{Z}.$$\end{theorem}
\section{Some Lemmas}
 \begin{lemma} {\rm(\cite{JLM09} Lemma 1)}
The coefficient $c(n,k)$ is an integer whose value only depends on
the congruence class of $k$ modulo $n$.\end{lemma}

Let $\kappa(m)=\prod_{p|m}p$ denote the squarefree kernel of $m$.

\begin{lemma} {\rm(\cite{JLM09} Corollary 1)}
We have $S(m)=S(\kappa(m))$ and $R(m)=R(\kappa(m))$.\end{lemma}

\begin{lemma}{\rm (Quantitative Form of Dirichlet¡¯s Theorem)}
Let $a$ and $m$ be coprime natural numbers and let $\pi(x;m,a)$
denote the number of primes $ p\le x$ that satisfy  $p\equiv
a\pmod{m}$. Then, as $x$ tends to infinity,
$$\pi(x; m, a)\sim\frac{x}{\varphi(m)\log x},$$where $\phi$ is Euler's toitent function.\end{lemma}

\begin{lemma} {\rm(\cite{JLM09} Corollary 2)}
Given $m,t\ge1$ and any real number $r>1$ , there exists a constant
$N_0(t,m,r)$ such that for every $n>N_0(t,m,r)$ the interval
$(n,rn)$ contains at least $t$ primes $p\equiv1\pmod{m}$.\end{lemma}

\section{ The proof of Theorem 1}

\begin{proof} We first prove that $S(m; s, t)=\mathbb{Z}$. Since
$S(m; s, t)= S(\kappa(m); s, t)$ and $S(m; s, t)\supseteq S(mp; s,
t)$, where $p\equiv1\pmod{s}$ is an odd prime, we may assume that
$m$ is square-free, $m>1$ and  $\mu(m)=1$. Suppose that $n>N_0(t,
ms,\frac{15}{8})$, then, by Lemma 2.4, there exist primes $p_1, p_2,
\ldots, p_t$ such that
$$N<p_1<p_2<\cdots<p_t<\frac{15}{8}n\quad \mbox{and} \quad
p_j\equiv1\pmod{ms}, \quad j=1, 2, \ldots, t.$$

Let $q_1, q_2$ be primes  such that $q_2>q_1>2p_1$, $q_1\equiv
t\pmod{s}$ and $q_2\equiv1\pmod{s}$ and put

\begin{equation} \label{eq:1}
m_1=\left\{ \begin{aligned}
         p_1p_2\cdots p_tq_1 & \quad \mbox{if} \, t\, \mbox{is even}; \\
                  p_1p_2\cdots p_tq_1q_2& \quad {\rm otherwise}.
                          \end{aligned} \right.
                          \end{equation}
Note that $m$ and $m_1$ are coprime, $m_1\equiv t\pmod{s}$ and that
$\mu(m_1)=-1$, where $\mu$ denotes the M\"obius function. Using
these observations we conclude that
\begin{equation} \label{eq:2}
\begin{split}
\Phi_{mm_1}(x)&\equiv\prod_{d|mm_1, d<2p_1}(1-x^d)^{\mu(\frac{mm_1}{d})}\pmod{x^{2p_1}} \\
 &\equiv\prod_{d|m}(1-x^d)^{\mu(\frac{m}{d})\mu(m_1)}
 \prod_{j=1}^t(1-x^{p_j})^{\mu(\frac{mm_1}{p_j})}\pmod{x^{2p_1}} \\
 &\equiv\Phi_m(x)^{\mu(m_1)}\prod_{j=1}^t(1-x^{p_j})^{-\mu(mm_1)}\pmod{x^{2p_1}} \\
 &\equiv\frac{1}{\Phi_m(x)}\prod_{j=1}^t(1-x^{p_j})^{\mu(m)}\pmod{x^{2p_1}} \\
 &\equiv\frac{1}{\Phi_m(x)}(1-\mu(m)(x^{p_1}+\cdots+x^{p_t}))\pmod{x^{2p_1}}.
 \end{split}
 \end{equation}
From (\ref{eq:2})  it follows that, if $p_t\le k<2p_1$, then
$$a(mm_1, k)=c(m, k)-\mu(m)\sum_{j=1}^tc(m, k-p_j).$$
By Lemma 2.1 we have $c(m, k-p_j)=c(m, k-1)$, and therefore
\begin{equation}\label{eq:3}
a(mm_1, k)=c(m, k)-\mu(m)tc(m, k-1)\,\, \mbox{with}\,\, p_t\le
k<2p_1.\end{equation} Since $\mu(m)=1$, we let $q_3<q_4$ be the
smallest two prime divisors of $m$. Here we also required that $n\ge
8q_4$, which ensures that $p_t+q_4<2p_1$. Note that
\begin{equation} \label{eq:4}
\begin{split}
\frac{1}{\Phi_{m}(x)}&\equiv\frac{(1-x^{q_3})(1-x^{q_4})}{1-x}\pmod{x^{q_4+2}} \\
  &\equiv1+x+x^2+\cdots+x^{q_3-1}-x^{q_4}-x^{q_4+1}\pmod{x^{q_4+2}}.
 \end{split}
 \end{equation}
 Thus $c(m, k)=1$ if $k\equiv\beta\pmod{m}$ with $\beta\in\{1, 2\}$
 and $c(m, k)=-1$ if $k\equiv\beta\pmod{m}$ with $\beta\in\{q_4,
 q_4+1\}$. This in combination with (\ref{eq:3}) shows that $a(m_1m,
 p_t+1)=1-t$ and $a(m_1m, p_t+q_4)=t-1$. Since $\{1-t,
 t-1|t\ge1\}=\mathbb{Z}$, then $S(m; s, t)=\mathbb{Z}$ and the first result follows.

 To prove $R(m; s, t)=\mathbb{Z}$. As before we may assume that $m>1$ is
 square-free and $\mu(m)=1$.

 Let $q_1, q_2$ be primes  such that $q_2>q_1>2p_1$, $q_1\equiv
t\pmod{s}$ and $q_2\equiv1\pmod{s}$ and put

\begin{equation} \label{eq:5}
\bar{m_1}=\left\{ \begin{aligned}
         p_1p_2\cdots p_tq_1q_2 & \quad \mbox{if} \, t\, \mbox{is even}; \\
                  p_1p_2\cdots p_tq_1& \quad {\rm otherwise}.
                          \end{aligned} \right.
                          \end{equation}
Note that $m$ and $m_1$ are coprime and that $\mu(\bar{m_1})=1$.
Reasoning as in the derivation of (\ref{eq:2}) we obtain
\begin{equation} \label{eq:6}
\frac{1}{\Phi_{mm_1}(x)}\equiv\frac{1}{\Phi_m(x)}(1-\mu(m)(x^{p_1}+\cdots+x^{p_t}))\pmod{x^{2p_1}}
 \end{equation} and from this $c(\bar{m_1}m, k)=a(m_1m, k)$ for
 $k\le 2p_1$. Reasoning as in the proof $S(m; s, t)=\mathbb{Z}$, we obtain $R(m; s, t)=\mathbb{Z}$. This completes the
 proof.
\end{proof}
{\bf Remark:} Since we do not need to consider the case $\mu(m)=-1$,
so a proof a little easier than that given in \cite{JLM09} is
obtained.

 {\bf Acknowledgments:} The  author is supported by NSF of China (No.
 10971072) and by
the Guangdong Provincial Natural Science Foundation (No.
8151027501000114).

\end{document}